\documentclass[11pt]{article}

\usepackage[active]{srcltx} 
\usepackage{amsmath,amsxtra,amssymb}
\usepackage{latexsym,array,verbatim,fancybox,multirow}
\usepackage{epsfig,graphicx,url,setspace}
\usepackage{algorithm,algorithmic}

\newtheorem{thm}{Theorem}

\newtheorem{prop}{Proposition}

\newtheorem{prob}{Problem}

\newcommand{\norm}[1]{\left\Vert#1\right\Vert}
\newcommand{\abs}[1]{\left\vert#1\right\vert}

\newcommand{\Real}{\mathbb R}

\newcommand{\mc}{\mathcal}

\begin{document}

\title{A Framework for Optimization under Limited Information}


\author{Tansu Alpcan \\
              Technical University Berlin \\
	      Deutsche Telekom Laboratories \\
              \textit{alpcan@sec.t-labs.tu-berlin.de}         
}

\date{}


\maketitle

\begin{abstract}
In many real world problems, optimization decisions have to be made with limited information. The decision maker may have no a priori or posteriori data about the often nonconvex objective function except from on a limited number of points that are obtained over time through costly observations. This paper presents an optimization framework that takes into account the information collection (observation), estimation (regression), and optimization (maximization) aspects in a holistic and structured manner.  Explicitly quantifying the information acquired at each optimization step using the entropy measure from information theory, the (nonconvex) objective function to be optimized (maximized) is modeled and estimated by adopting a Bayesian approach and using Gaussian processes as a state-of-the-art regression method. The resulting iterative scheme allows the decision maker to solve the problem by expressing preferences for each aspect quantitatively and concurrently. 

\end{abstract}

\section{Introduction} \label{sec:intro}

In many real world problems, optimization decisions have to be made with limited information. Whether it is a static optimization or dynamic control problem, obtaining detailed and accurate information about the problem or system can often be a costly and time consuming process. In some cases, acquiring extensive information on system characteristics may be simply infeasible. In others, the observed system may be so nonstationary that by the time the information is obtained, it is already outdated due to system's fast-changing nature. Therefore, the only option left to the decision-maker is to develop a strategy for collecting information efficiently and choose a model to estimate the ``missing portions'' of the problem in order to solve it satisfactorily and according to a given objective.


To make the discussion more concrete, consider the problem of maximizing a (Lipschitz) continuous \textit{nonconvex} objective function, which is unknown except from its value at only a small number of data points. The decision maker may have no a priori information about the function and start with zero data points. Furthermore, only a limited number of --possibly noisy-- observations may be available before making a decision on the maximum value and its location. The function itself, however, remains unknown even after the decision is made. \textit{What is the best strategy to address this problem}? 

The decision making framework presented in this paper captures the posed problem by taking into account the information collection (observation), estimation (regression), and (multi-objective) optimization aspects in a holistic and structured manner. Hence, the framework enables the decision maker to solve the problem by expressing preferences for each aspect quantitatively and concurrently. It explicitly incorporates many concepts that have been implicitly considered by heuristic schemes, and  builds upon many results from seemingly disjoint but relevant fields such as information theory, machine learning, and optimization and control theories.  Specifically, it combines concepts from these fields by
\begin{itemize}
 \item explicitly quantifying the information acquired using the entropy measure from information theory,
 \item modeling and estimating the (nonconvex) function or (nonlinear) system adopting a Bayesian approach and using Gaussian processes as a state-of-the-art regression method, 
 \item using an iterative scheme for observation, learning, and optimization,
 \item capturing all of these aspects under the umbrella of a multi-objective ``meta'' optimization formulation.
\end{itemize}

Despite methods and approaches from machine (statistical) learning are heavily utilized in this framework, the problem at hand is very different from many classical machine learning ones, even in its learning aspect. In most classical application domains of
machine learning such as data mining, computer vision, or image and voice recognition, the difficulty is often in handling significant amount of data in contrast to lack of it. Many methods such as Expectation-Maximization (EM) inherently make this assumption, except from ``active learning'' schemes \cite{Bishopbook}. Information
theory plays plays an important role in evaluating scarce (and expensive) data and developing strategies for obtaining it. Interestingly, data scarcity converts at the same time the disadvantages of some methods into advantages, e.g. the scalability problem of Gaussian processes.

It is worth noting that the class of problems described here are much more frequently encountered in practice than it may first seem. For example, the class of black-box methods known as ``kriging'' \cite{kriging1} have been applied to such problems in geology and mining as well as to hydrology since mid-1960s. In addition, the solution framework proposed is applicable to a wide variety of fields due to its fundamental nature. One example is decentralized resource allocation decisions in networked and complex systems, e.g. wired and wireless networks, where parameters change quickly and global information on network characteristics are not available at the local decision-making nodes. Another example is security-related decisions where opponents spend a conscious effort to hide their actions. A related area is security and information technology risk management in large-scale organizations, where acquiring information on individual subsystems and processes can be very costly. Yet another example application is in biological systems where individual organisms or subsystems operate autonomously (even if they are part of a larger system) under limited local information.

\section{Problem Definition and Approach} \label{sec:problem}

A concrete definition of the motivating problem mentioned in the introduction section is helpful for describing the multiple aspects of the limited information decision making framework. Without loss of any generality, let 
$$\mc X  \subseteq \Psi  \subset \Real^{d}$$ 
be a nonempty, convex, and compact (closed and bounded) subset of the original problem domain $\Psi$ of $d$ dimensions. The original domain $\Psi$ does not have to be convex, compact, or even fully known. However, adopting a ``divide and conquer'' approach, the subset $\mc X$ provides a reasonable starting point. 
Define next the objective function to be maximized 
$$f: \mc X \rightarrow \Real, $$
which is unknown except from on a finite number of points (possibly imperfectly) observed. As a simplifying
assumption, let $f$ be Lipschitz continuous on $\mc X$. One of the main distinguishing characteristics of this problem is the limitations on set of observations 
$$\Omega_n:=\{x_1,\ldots,x_n \, : x_i \in \mc X\; \forall i,\; n \geq 1 \},$$ 
due to cost of obtaining information or non-stationarity of the underlying system. Assume for now that the cost of observing the value of the objective function $f(x)$ is the same for any $x \in \mc X$. Then, a basic search problem is defined as follows:

\begin{prob}[\textit{Basic Search Problem}] \label{prob:search1}
Consider a Lipschitz-continuous objective function $f: \mc X \rightarrow \Real$ on the $d$-dimensional nonempty, convex, and compact set $\mc X \subset \Real^{d}$. The function is unknown except from on a finite number of observed data points. What is the best search  strategy 
$$\Omega_N:=\{x_1,\ldots,x_N \, : x_i \in \mc X\; \forall i,\; N \geq 1 \}$$ 
that solves
$$ \max_{ \Omega_N} f(x) ,$$
for a given $N$?
\end{prob}

The number of observations, $N$, in Problem~\ref{prob:search1} may be imposed by the nature of the specific application domain. In many problems, where there is no time constraint, adopting an iterative (one-by-one) approach, and hence choosing $N=1$ is clearly beneficial 
as it allows for usage of incoming new information at each step. Alternatively, the assumption on the equal observation cost can be relaxed and be formulated as a constraint
$$ \sum_{x \in \Omega_n} c_o(x) \leq C, $$
where $c_o(x): \mc X \rightarrow \Real $ is the observation cost function, and the scalar $C$ is the total ``exploration budget''. It is also possible to define this cost iteratively based on the (distance from) previous observation, e.g. $c_o(x_n,x_{n-1})$. In such cases, a location-based iterative search scheme can be considered.

The simplest (both conceptually and computationally) strategy to solve Problem~\ref{prob:search1} is random search on the domain $\mc X$. As such no attempt is made to ``learn'' the properties of the function $f$. Unless, $f$ is 
``algorithmically random'' \cite{algobook}, which is rarely the case, this strategy wastes the information collected on $f$. A slightly more complicated and very popular set of strategies combine random search with simple modeling of the function through gradient methods. In this case, the collected information is used to model  $f$ rudimentarily using derived gradients to ``define slopes'' in a heuristic manner. Then, these slopes of $f$ are explored step-by-step in the upwards direction to find a local maximum, after which the search algorithm randomly jumps to another location. It is also possible to randomize the gradient climbing scheme for additional flexibility \cite{simannealing}.


The framework presented in this paper takes one further step and \textbf{explicitly} models the (entire) objective function $f$ (on the set $\mc X$) using the information collected instead of heuristically describing only the slopes. The function $\hat f$, which models, approximates, and estimates $f$, belongs to a certain class functions such that $\hat f \in \mc F$. The selection and properties of this class is based on ``a priori'' information available and can be interpreted as the ``world view'' of the decision maker. These properties can often be expressed using meta-parameters which are then updated based on the observations through a separate optimization process. Likewise, a slower time-scale process can be used for model selection if processing capabilities permit a multi-model approach.

This model-based search process, which lies at the center of the framework, is fundamentally a manifestation of the Bayesian approach \cite{MacKaybook}. It first imposes explicit and a priori modeling assumptions by choosing $\hat f$ from a certain class of functions, $\mc F$, and then infers (learns, updates) $\hat f$ in a structured manner as more information becomes available through observations. 

From a computational point of view, the decision making framework with limited information lies at one end of the computation vs. observation spectrum, while random search is at the opposite end. The framework tries to utilize each piece of information to the maximum possible extent almost regardless of the computational cost. The underlying assumption here is: \textbf{observation is very costly whereas computation is rather cheap}. This assumption is not only valid for a wide variety of problems from different fields ranging from networking and security to economics and risk management, but also inspired from biological systems. In many biological organisms, from single cells to human beings, operating close to this end of the computation-observation spectrum is more advantageous than doing random search. 

When doing random search on the domain $\mc X$, at each stage i.e. given the previous observations, each remaining candidate data point provides equivalent amount of information. However, this is not the case when doing model-based search. Depending on the model adopted and previous information collected, different unexplored points provide different amount of information. This information can be exactly quantified using the definition of entropy and information from the field of (Shannon) information theory. Accordingly, the scalar quantity $\mc I(\hat f, \Omega_n)$ denotes the aggregate information obtained from the set of observations $\Omega_n$ within the model represented by $\hat f$. A related issue is the reliability and possibly noisy nature of observations, which will be discussed in further detail in the next section.

An extension of Problem~\ref{prob:search1} that captures the aspects discussed above is defined next.
\begin{prob}[\textit{Model-based Search Problem}] \label{prob:search2}
Let $f: \mc X \rightarrow \Real$ be an objective function on the $d$-dimensional nonempty, convex, and compact set $\mc X \subset \Real^{d}$, which is unknown except from on a finite number of observed data points. Further let $\hat f(x)$ be an estimate of the objective function obtained using an a priori model and observed data. What is the best search  strategy $\Omega_N:=\{x_1,\ldots,x_N \, : x_i \in \mc X\; \forall i,\; N \geq 1 \}$ that solves the multi-objective problem with the following components?
\begin{itemize}
 \item \textit{Objective 1:} $ \max_{\Omega_N} f(x) \text{ given } \hat f(x)$
 \item \textit{Objective 2:} $ \arg \min_{\Omega_N} R\left( f(x), \hat f(x) \right) , \; \hat f \in \mc F$
 \item \textit{Objective 3:} $ \max_{\Omega_N} \mc I(\hat f, \Omega_n)$
\end{itemize}
Here, $R(\cdot,\cdot)$ is a risk or expected loss function  quantifying the mismatch between actual and estimated functions on the observation data \cite{GPbook}. The scalar quantity $\mc I$ is the aggregate information obtained from the set of observations $\Omega_N$ within the model represented by $\hat f$. The cardinality of $\Omega_N$, $N$, can be either given, e.g. $N=1$, or defined as an additional constraint $ \sum_{x \in \Omega_n} c_o(x) \leq C$, where $c_o(x): \mc X \rightarrow \Real $ is the observation cost function, and the scalar $C$ is the total ``exploration budget''.
\end{prob}


\begin{figure}[htp]
 \centering
 \includegraphics[width=0.6\columnwidth]{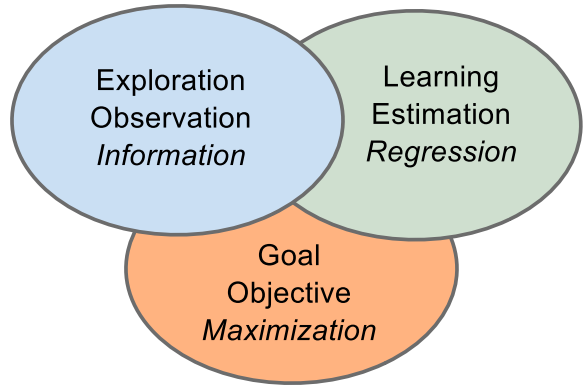}
 \caption{The three fundamental aspects of decision making with limited information.}
 \label{fig:goals1}
\end{figure}

It is important to observe here that the three objectives defined in Problem~\ref{prob:search2} are (almost) independent from and orthogonal to each other despite being closely related. \textit{Objective 1} purely aims to maximize the unknown objective function $f$ using the best estimate (model) $\hat f$. \textit{Objective 2} focuses on minimizing the error between the estimate $\hat f$ and the real unknown function $f$ based on the observations made. \textit{Objective 3} tries to maximize the amount of information provided by each (costly) observation or experiment. It is worth noting that 
\textit{Objective 3} is independently formulated from \textit{Objective 2}, in other words, exploration is done independently from estimation. In contrast, ensuring a balance between \textit{Objective 1} and \textit{2} is necessary to ensure that solution is robust. These objectives and the fundamental aspects of decision making with limited information are visually depicted in Figure~\ref{fig:goals1}. 

\begin{table}[htp]
  \caption{Fundamental Trade-offs}
\begin{center}
    \begin{tabular}{l|c|l}
     \hline \hline  
     Exploration &  & Exploitation \\[1 ex]
     Observation & versus & Computation \\[1 ex]
     Robustness  &  & Optimization \\ [1 ex]
     \hline \hline
    \end{tabular}
  \label{tbl:tradeoff}
\end{center}
 \end{table}
There are multiple trade-offs that are inherent to this problem as listed in Table~\ref{tbl:tradeoff}. The first one, exploration versus exploitation, puts exploration or obtaining more observations against exploitation, i.e. trying to achieve the given objective. Observation versus computation captures the trade-off between building sophisticated models using the available information to the fullest extend and making more observations. Robustness versus optimization puts risk avoidance against optimization with respect to the original objective as in exploitation.

\section{Methodology} \label{sec:methods}

This section presents the methods that are utilized within the framework which addresses the problem defined in the previous one. First, the regression model and Gaussian Processes (GP) are presented. Subsequently, modeling and measurement of information is discussed based on (Shannon) information theory. 

\subsection{Regression and Gaussian Processes (GP)} \label{sec:gp}

Problem~\ref{prob:search2} presented in the previous section involves inferring or learning the function $f$ using the set of observed data points. This is known as the \textit{regression} problem in machine learning and is a supervised learning method since the observed data constitutes at the same time the learning data set. This learning process involves selection of a ``model'', where the learned function $\hat f$ is, for example, expressed in terms of a set of parameters and specific basis functions, and at the same time minimization of an error measure between the functions $f$ and $\hat f$ on the learning data set. Gaussian processes (GP) provide a nonparametric alternative to this but follow in spirit the same idea. 

The main goal of regression involves a trade-off. On the one hand, it tries to minimize the \textit{observed} error between $f$ and $\hat f$. On the other, it tries to infer the ``real'' shape of $f$ and make good estimations using 
$\hat f$ even at unobserved points. If the former is overly emphasized, then one ends up with ``over fitting'', which means $\hat f$ follows $f$ closely at observed points but has weak predictive value at unobserved ones. This delicate balance is usually achieved by balancing the prior ``beliefs'' on the nature of the function, captured by the model (basis functions), and fitting the model to the observed data. 

This paper focuses on Gaussian Process \cite{GPbook} as the chosen regression method within the framework developed without loss of any generality. There are multiple reasons behind this preference. Firstly, GP provides an elegant mathematical method for easily combining many aspects of the framework. Secondly, being a nonparametric method GP eliminates any discussion on model degree. Thirdly, it is easy to implement and understand as it is based on well-known Gaussian probability concepts. Fourthly, noise in observations is immediately taken into account if it is modeled as Gaussian. Finally, one of the main drawbacks of GP namely being computational heavy, does not really apply to the problem at hand since the amount of data available is already very limited.

It is not possible to present here a comprehensive treatment of GP. Therefore, a very rudimentary overview is provided next within the context of the decision making problem. Consider a set of $M$ data points 
$$\mc D=\{x_1, \ldots, x_M\},$$
where each $x_i \in \mc X$ is a $d-$dimensional vector, and the corresponding vector of scalar values is $f(x_i), \; i=1,\ldots,M$. Assume that the observations are distorted by a zero-mean Gaussian noise, $n$ with variance $\sigma \sim \mc N(0,\sigma)$. Then, the resulting observations is a vector of Gaussian $y=f(x)+n \sim \mc N(f(x),\sigma)$. 

A GP is formally defined as a collection of random variables, any finite number of which have a joint Gaussian distribution. It is completely specified by its mean function $m(x)$ and covariance function $C(x,\tilde x)$, where
$$ m(x)=E[\hat f(x)] \text{ and } C(x,\tilde x)=E[(\hat f(x)-m(x))(\hat f(\tilde x)-m(\tilde x))], 
\; \forall x, \tilde x \in \mc D. $$

Let us for simplicity choose $m(x)=0$. Then, the GP is characterized entirely by its covariance function $C(x,\tilde x)$. Since the noise in observation vector $y$ is also Gaussian, the covariance function can be defined as the sum of a \textit{kernel function} $Q (x,\tilde x)$ and the diagonal noise variance 
\begin{equation} \label{e:gcov}
 C(x,\tilde x) = Q (x,\tilde x) + \sigma I, \; \forall \, x, \tilde x \in \mc D ,
\end{equation}
where $I$ is the identity matrix. While it is possible to choose here any (positive definite) kernel $Q(\cdot,\cdot)$, one classical choice is 
\begin{equation} \label{e:gaussiankernel}
 Q(x,\tilde x)=\exp \left[-\frac{1}{2}\norm{x -\tilde x}^2 \right].
\end{equation}
Note that GP makes use of the well-known \textit{kernel trick} here by representing an infinite dimensional continuous function
using a (finite) set of continuous basis functions and associated vector of real parameters in 
accordance with the \textit{representer theorem} \cite{schoelkopfbook}.

The (noisy)\footnote{The special case of perfect observation without noise is handled the same way as long as the kernel function $Q(\cdot,\cdot)$ is positive definite} training set $(\mc D, y)$ is used to define the corresponding GP, $\mc{GP} (0,C(\mc D))$, through the $M \times M$ covariance function $C(\mc D)=Q+\sigma I$, where the conditional Gaussian distribution of any point outside the training set, $\bar y \in \mc X, \bar y \notin \mc D$, given the training data $(\mc D, t)$  can be computed as follows. Define the vector 
\begin{equation} \label{e:k}
 k(\bar x)=[Q(x_1,\bar x), \ldots Q(x_M,\bar x)]
\end{equation}
and scalar 
\begin{equation} \label{e:kappa}
\kappa=Q(\bar x,\bar x)+\sigma.
\end{equation}
Then, the conditional distribution  $p(\bar y | y)$ that characterizes the $\mc{GP} (0,C)$ is a Gaussian $\mc N(\hat f,v)$ with mean $\hat f$ and variance $v$,
\begin{equation} \label{e:gp1}
 \hat f(\bar x)=k^T C^{-1} y \text{ and } v(\bar x)=\kappa - k^T C^{-1} k .
\end{equation}

This is a key result that defines GP regression as the mean function $\hat f(x)$ of the Gaussian distribution and provides a prediction of the objective function $f(x)$. At the same time, it belongs to the well-defined class $\hat f \in\mc F$, which is the set of all possible sample functions of the GP
$$\mc F := \{\hat f(x): \mc X \rightarrow \Real \text{ such that }  \hat f \in \mc{GP} (0,C(\mc D)),\; \forall \mc D, \, C \} ,$$
where $ C(\mc D)$ is defined in (\ref{e:gcov}) and $\mc{GP}$ through (\ref{e:k}), (\ref{e:kappa}), and (\ref{e:gp1}), above.
Furthermore, the variance function $v(x)$ can be used to measure the uncertainty level of the predictions provided by $\hat f$, which will be discussed in the next subsection.


\subsection{Quantifying Information in Observations} \label{sec:obsinfo}

In the framework presented, each observation provides a data point to the regression problem (estimating $f$ by constructing $\hat f$) as discussed in the previous subsection. Many works in the learning literature consider the ``training'' data used in regression available (all at once or sequentially) and do not discuss the possibility of the decision maker influencing or even optimizing the data collection process. The \textit{active learning} problem defined in Section~\ref{sec:problem} requires, however, exactly addressing the question of ``how to quantify information obtained and optimize the observation process?''. Following
the approach discussed in \cite{MacKaydataselect,MacKaybook}, the framework here provides a precise answer to this question.

Making any decision on the next (set of) observations in a principled manner necessitates first \textit{measuring the information obtained from each observation within the adopted model}. It is important to note that the information measure here is dependent on the chosen model. For example, the same observation provides a different amount of information to a random search model than a GP one.

Shannon information theory readily provides the necessary mathematical framework for measuring the information content of a variable. Let $p$ be a probability distribution over the set of possible values of a discrete random variable $A$. The \textbf{entropy} of the random variable is given by
$H(A)=\sum_i p_i \log_2 (1/p_i)$, which quantifies the amount of uncertainty. Then, the information obtained from an observation on the variable, i.e. reduction in uncertainty, can be quantified simply by taking the difference of its initial and final entropy, 
$$\mc I=H_0 - H_1. $$ 
It is important here to avoid the common conceptual pitfall of equating entropy to information itself as it is sometimes done in communication theory literature.\footnote{Since this issue is not of great importance for the class of problems considered in communication theory, it is often ignored. However, the difference is of conceptual importance in this problem. See
\url{http://www.ccrnp.ncifcrf.gov/~toms/information.is.not.uncertainty.html} for a detailed discussion.} 
Within this framework, (Shannon) \textit{information is defined as a measure of the decrease of uncertainty after (each) observation (within a given model)}. This can be best explained with the following simple example.

\subsubsection{Example: Bisection} \label{exp:bisection}
Choose a number between $1$ and $64$ randomly with uniform probability (prior). What is the best searching strategy for finding this number? Let the random variable $A$ represent this number. In the beginning the entropy of $A$ is 
$$H_0(A)=\sum_{i=1}^{64} \dfrac{1}{64} \log_2 \left( \dfrac{1}{64} \right)=6 \text{ (bits)}.$$ 
The information maximization problem is defined as
$$ \max \mc I= \max H_0 - H_1 = \min H_1 ,$$
since $H_0$, the entropy before the action (obtaining information) is constant. The entropy $H_1$ is the one after information is obtained, and hence is directly affected by the specific action chosen. Now, define the action as setting a threshold $1 < t <64$ to check whether the chosen number is less or higher than this threshold $t$. To simplify the analysis, consider a continuous version of the problem by defining $p$ as the probability of the chosen number being less than the threshold. Thus, in this uniform prior case, the problem simplifies to
$$ \min_p H_1= \min_p \; p \log(p) + (1-p) \log (1-p),$$
which has the derivative 
$$ \dfrac{d H_1}{d p}= \log(p) - \log(1-p) .$$
Clearly, the threshold $p^*=0.5$ is the global minimum, which roughly corresponds to $t=32$ (ignoring quantization and boundary effects). Thus, bisection from the middle is the optimal search strategy for the uniform prior. In this example, the number can be found in the worst-case in $6$ steps, each providing one bit of information. Nonuniform probabilities (priors) can be handled in a similar way. 

If this search process (bisection) is repeatedly applied without any feedback, then it results in the optimal quantization of the search space both in the uniform case above and for the nonuniform probabilities. If feedback is available, i.e. one learns after each bisection whether the number is larger or less than the boundary, then this is as shown the best search strategy.

\section{Model} \label{sec:model}

\begin{figure}[htp]
 \centering
 \includegraphics[width=0.9\columnwidth]{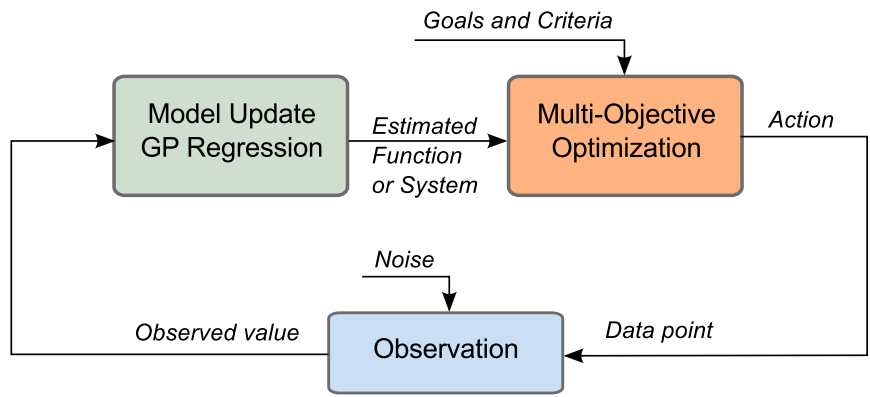}
 \caption{The main parts of the underlying model of the decision making framework.}
 \label{fig:model1}
\end{figure}

The model adopted in the framework for decision making with limited information builds on the methods presented in the previous section and addresses the problem introduced in Section~\ref{sec:problem}. The model consists of three main parts: observation, update of GP for regression, and optimization to determine next action. These three steps, shown in Figure~\ref{fig:model1} are taken iteratively to achieve the objectives in Problem~\ref{prob:search2}. As a result of its iterative nature, this approach can be considered in a sense similar to the well-known Expectation-Maximization algorithm \cite{Bishopbook}.

Observations, given that they are a scarce resource in the class of problems considered, play an important role in the model. Uncertainties in the observed quantities can be modeled as additive noise. Likewise, properties (variance or bias) of additive noise can be used to model the reliability of (and bias in) the data points observed. GPs provide a straightforward mathematical structure for incorporating these aspects to the model under some simplifying assumptions.

The set of observations collected provide the (supervised) training data for GP regression in order to estimate the characteristics of the function or system at hand. This process relies on the GP methods described in Subsection~\ref{sec:gp}. Thus, at each iteration an up-to-date description of the function or system is obtained
based on the latest observations. Specifically, $\hat f$ provides an estimate of the original function $f$.\footnote{See \cite[Chap 7.2]{GPbook} for a discussion on asymptotic analysis of GP regression. It should not be noted, however, that asymptotic properties are of little relevance to the problem at hand. }
Assuming an additive Gaussian noise model, the noise variance $\sigma$ can be used to model uncertainties, e.g. older and noisy data resulting in higher $\sigma$ values.

The final and most important part of the model provides a basis for determining the next action after an optimization process that takes into account all three objectives in Problem~\ref{prob:search2}. The information aspect of these objectives is already discussed in Subsection~\ref{sec:obsinfo}. An important issue here is the fact that there are infinitely many candidate points in this optimization process, but in practice only a finite collection of them can be evaluated. 

\subsection{Sampling Solution Candidates}

When making a decision on the next action through multi-objective optimization, there are (infinitely) many candidate points. A pragmatic solution to the problem of finding solution candidates is to (adaptively) sample the problem domain $\mc X$ to obtain the set  
$$\Theta:=\{x_1, \ldots, x_T : x_i \in \mc X, \, x_i \notin \mc D, \; \forall i \}$$ 
that does not overlap with known points. In low (one or two) dimensions, this can be easily achieved through grid
sampling methods. In higher dimensions, (Quasi) Monte Carlo schemes can be utilized. For large problem domains, the current domain of interest $\mc X$ can be defined around the last or most promising observation in such a way that such a sampling is computationally feasible. 
Likewise, multi-resolution schemes can also be deployed to increase computational efficiency.

Although such a solution may seem restrictive at first glance, it is in spirit not very different from other schemes such as simulated annealing, which are widely used to address nonconvex optimization problems. However, a major difference
between this and other schemes is the fact that the candidate sampling and evaluation are done here ``a priori'' due to experimentation being costly while other methods rely on abundance of information.

A natural question that arises is: whether and under what conditions does such a sampling method give satisfactory results. 
The following result from \cite{tempo-sampling,tempobook} provides an answer to this question in terms of number of samples required.
\begin{thm} \label{thm:sampling}
Define a multivariate function $f(x)$ on the convex, compact set $\mc X$, which admits the maximum $x^*=\arg \max_{x \in \mc X} f(x)$. Based on a set of $N$ random samples $\Theta=\{x_1, \ldots, x_N: x_i \in \mc X \; \forall i \}$ from the entire set $\mc X$, let $\hat x:= \arg \max_{x \in \Theta} f(x)$ be an estimate of the maximum $x^*$. 

Given an $\varepsilon>0$ and $\delta>0$, the minimum number of random samples $N$ which guarantees that
$$ Pr\left(  Pr[f(x^*)>f(\hat x)] \leq \varepsilon \right) \geq 1-\delta,$$
i.e. the probability that 'the probability of the real maximum surpassing the estimated one being less than $\varepsilon$' is larger than $1-\delta$, is
$$ N \geq \dfrac{\ln 1/ \delta}{ 1/ (1-\varepsilon)} .$$
Furthermore, this bound is tight if the function $f$ is continuous on $\mc X$.
\end{thm}
It is interesting and important to note that this bound is independent of the sampling distribution used (as long as it covers the whole set $\mc X$ with nonzero probability), the function $f$ itself, as well as the properties and dimension of the set $\mc X$.
 
\subsection{Quantifying Information in GP}

The information measurement and GP approaches in Section~\ref{sec:methods} can be directly combined. Let the zero-mean multivariate Gaussian (normal) probability distribution be denoted as
\begin{equation} \label{e:multivargauss}
 p(x)=\dfrac{1}{\sqrt{2\pi |C_p(x)}|} \exp \left( -\frac{1}{2}[x-m]^T|C_p(x)|^{-1} [x-m]\right) ,\;  x \in \mc X,
\end{equation}
where $|\cdot|$ is the determinant, $m$ is the mean (vector) as defined in (\ref{e:gp1}), and $C_p(x)$ is the covariance matrix as a function of the newly observed point $x \in \mc X$ given by
\begin{equation} \label{e:covx}
 C_p(x)=\left[ 
\begin{array}{cccc}
 &   &  &  \\
 & C(\mc D) &  & k(x)^T \\
 &   &  &  \\
 & k(x) &  & \kappa
\end{array}
\right] .
\end{equation}
Here, the vector $k(x)$ is defined in (\ref{e:k}) and $\kappa$ in (\ref{e:kappa}), respectively. The matrix $C(\mc D)$ is the covariance matrix based on the training data $\mc D$ as defined in (\ref{e:gcov}). 

The entropy of the multivariate Gaussian distribution (\ref{e:multivargauss}) is \cite{entropygaussian} 
$$ H(x)=\dfrac{d}{2}+\dfrac{d}{2}\ln(2\pi)+\dfrac{1}{2} \ln |C_p(x)| ,$$
where $d$ is the dimension. Note that, this is the entropy of the GP estimate at the point $x$ based on the available data $\mc D$. The aggregate entropy of the function on the region $\mc X$ is given by
\begin{equation} \label{e:aggentropy}
 H^{agg}:=\int_{x \in \mc X} \dfrac{1}{2} \ln |C_p(x)| dx.
\end{equation}

The problem of choosing  a new data point $\hat x$ such that the information obtained from it within the
GP regression model is maximized can be formulated in a way similar to the one in the bisection example:
\begin{equation} \label{e:infocollect1}
 \hat x=\arg \max_{\tilde x} \mc I= \arg \max_{\tilde x} \int_{x \in \mc X} \left[ H_0 - H_1 \right] \, dx = \arg \min_{\tilde x}  \int_{x \in \mc X} \dfrac{1}{2} \ln |C_q(x,\tilde x)| dx,
\end{equation}
where the integral is computed over all $x \in \mc X$, and the covariance matrix $C_q(x, \tilde x)$ is defined as
\begin{equation} \label{e:covxbar}
 C_q(x, \tilde x)=\left[ 
\begin{array}{ccccc}
 &   &         &          &     \\
 & C(\mc D)&   & k^T(\tilde x) & k^T(x) \\
 &   &         &          &  \\
 &  k(\tilde x) &   & \tilde \kappa & Q(x,\tilde x) \\
 &  k(x) &  & Q(x,\tilde x) & \kappa
\end{array}
\right] ,
\end{equation}
and $\tilde \kappa=Q(\tilde x,\tilde x)+\sigma$. Here, $C(\mc D)$ is a $M \times M$ matrix and $C_q$ is a $(M+2) \times (M+2)$ one, whereas $\kappa$ and  $Q(x,\tilde x)$ are scalars and $k$ is a $M \times 1$ vector. 
This result is summarized in the following proposition.

\begin{prop} \label{thm:GPsearch}
As a maximum information data collection strategy for a Gaussian Process with a covariance matrix $C(\mc D)$, the next observation $\hat x$ should be chosen in such a way that 
$$ \hat x= \arg \max_{\tilde x} \mc I= \arg  \min_{\tilde x}  \int_{x \in \mc X} \ln |C_q(x,\tilde x)| dx,$$
where $C_q(x, \tilde x)$ is defined in (\ref{e:covxbar}).
\end{prop}


\subsubsection*{An Approximate Solution to Information Maximization}

Given a set of (candidate) points $\Theta$ sampled from $\mc X$, the result in Proposition~\ref{thm:GPsearch} can be revisited. The problem in (\ref{e:infocollect1}) is then approximated  \cite{tempobook} by
\begin{eqnarray} \label{e:infocollect2}
 \max_{\tilde x} \mc I \approx \min_{\tilde x} \sum_{x \in \Theta} \ln |C_q(x,\tilde x)| \\
 \nonumber \Rightarrow  \hat x= \arg \min_{\tilde x \in \Theta}  \prod_{x \in \Theta} |C_q(x, \tilde x)|,
\end{eqnarray}
using monotonicity property of the natural logarithm and the fact that the determinant of a covariance matrix is non-negative.
Thus, the following counterpart of Proposition~\ref{thm:GPsearch}
is obtained:
\begin{prop} \label{thm:GPsearch2}
As an approximately maximum information data collection strategy for a Gaussian Process with a covariance matrix $C(\mc D)$ and given a collection of candidate points $\Theta$, the next observation $\hat x \in \Theta$ should be chosen in such a way that 
$$ \hat x=  \arg  \min_{\tilde x \in \Theta}  \prod_{x \in \Theta} |C_q(x, \tilde x)| \approx \arg \max_{\tilde x \in \Theta} \mc I,$$
where $C_q(x, \tilde x)$ is given in (\ref{e:covxbar}).
\end{prop}

Although it is an approximation, finding a solution to the optimization problem in Proposition~\ref{thm:GPsearch2} can still be computationally  costly. Therefore, a greedy algorithm is proposed as a computationally simpler alternative. 
Let,  $x^* \in \Theta$ be defined as 
$$ x^* :=\arg \max_{x \in \Theta} |C_p(x)|=|C(\mc D)|\, |\kappa(x) - k(x) C^{-1}(\mc D) k^T(x) |,$$
where the matrix $C_p$ is given by (\ref{e:covx}) \cite{matrixcookbook}.  The first term above, $|C(\mc D)|$ is fixed and the second one, 
$$|\kappa(x) - k(x) C^{-1}(\mc D) k^T(x) |, $$
is the same as the GP variance $v(x)$ in (\ref{e:gp1}). Hence, the sample $x^*$ is one of those with the maximum variance in the set $\Theta$, given current data $\mc D$. 

It follows from (\ref{e:covxbar}) and basic matrix theory that if $\tilde x=x$ for a given $x$ then
$ |C_q(x, \tilde x)|$ is minimized. As a simplification, ignore the dependencies between $C_q(x, \tilde x)$ matrices for different $x \in \Theta$. Then, choosing the maximum variance $\hat x$ as
$$ \hat x = \arg \max_{\tilde x \in \Theta} v(\tilde x) \approx \arg \min_{\tilde x \in \Theta}  \prod_{x \in \Theta} |C_q(x, \tilde x)|,$$
leads to a large (possibly largest) reduction in $\prod_{x \in \Theta} |C_q(x, \hat x)|$, and hence
provides a rough approximate solution to (\ref{e:infocollect2}) and to the result in Proposition~\ref{thm:GPsearch}.
This result is consistent with widely-known heuristics such as ``maximum entropy'' or ``minimum variance'' methods \cite{activelearning} and a variant has been discussed in \cite{MacKaydataselect}.

\begin{prop} \label{thm:GPsearch3}
Given a Gaussian Process with a covariance matrix $C(\mc D)$ and a collection of candidate points $\Theta$, an approximate solution to the maximum information data collection problem defined in Proposition~\ref{thm:GPsearch} is to choose
the sample point(s) $\tilde x$ in such a way that it has (they have) the maximum variance within the set $\Theta$.
\end{prop}


%

\section{Optimization with Limited Information} \label{sec:staticopt}

%


Let $f: \mc X \rightarrow \Real$ be the unknown Lipschitz-continuous function of interest on the $d$-dimensional nonempty, convex, and compact set $\mc X \subset \Real^{d}$. The amount of information about this function available to the decision maker is limited to a finite number of possibly noisy observations. Since the observations are costly, the goal of the decision maker is to find the maximum of $f$, estimate $f$ as accurately as possible using available observations, and select the most informative data points, at the same time. This naturally calls for an iterative and myopic optimization procedure since each new observation provides a new data point that concurrently affects the maximization, function estimation (regression), and
information quantity.

The first and basic objective is the maximization of the function $f(x)$ on $x \in \mc X$. As a simplification, observations are assumed to be sequential, one at a time. Since $f$ is basically unknown, this problem has to be formulated as 
$$ \max_{\tilde x \in \mc X} F_1(\tilde x)= \hat f(\tilde x),$$
where $\hat f$ is the best estimate obtained through GP regression (\ref{e:gp1}) using the current data set $\mc D$. Data uncertainty (observation errors) is modeled through additive Gaussian noise with variance $\sigma$ as a first approximation. 

The second objective is to minimize the difference (estimation error) between $\hat f$ and  $f$. Define $e(x)=\hat f(x) - f(x), \; \forall x \in \mc X$. Given the set of noisy observations 
$$\mc O=\{f(x_i)+n(x_i): x \in \mc D, \, \forall i \} ,$$
where $n \sim \mc N(0,\sigma)$ denotes zero mean Gaussian noise, it is possible to use another GP regression (\ref{e:gp1}) to estimate this error function, $\hat e(\mc D, x)$, on the entire set $\mc X$. Thus, the second objective is to ensure that the next observation $\tilde x$ solves
$$ \min_{\tilde x \in \mc X} F_2(\tilde x)= \int_{\tau \in \mc X} \abs{\hat e(\tilde x,\mc D, \tau)} d \tau.$$
Note that, $F_2$ here corresponds to a risk or loss estimate function.

The third objective is to maximize the amount of information obtained with each observation $\tilde x$, or
$$ \max_{\tilde x \in \mc X} F_3(\tilde x)=\mc I (\tilde x, \hat f)=  \int_{x \in \mc X} \ln |C_q(x,\tilde x)| dx, $$
given the best estimate of the original function, $\hat f$. This objective has already been discussed in Section~\ref{sec:obsinfo} in detail.

\begin{figure}[htp]
 \centering
 \includegraphics[width=\columnwidth]{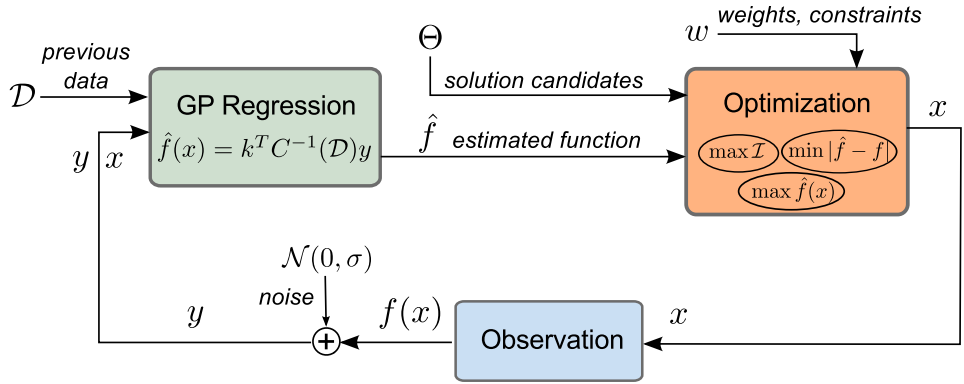}
 \caption{The decision making framework for static optimization with limited information.}
 \label{fig:model2}
\end{figure}

The values of the three objectives, $F_1,\, F_2,\, F_3$, cannot be evaluated numerically on the entire set $\mc X$. Therefore, a sampling method is used as described in Section~\ref{sec:model} to obtain a set of solution candidates $\Theta$, which replaces $\mc X$ in the maximization and minimization problems above. Next, specific problem formulations are presented based on such a sampling of solution candidates. The overall structure of the framework is visualized in Figure~\ref{fig:model2}. 

\subsection{Solution Approaches}

The most common approach to multi-objective optimization is the \textbf{weighted sum method}  \cite{moosurvey1,moosurvey2}. The three objectives discussed above can be combined to obtain a single objective using the respective weights $[w_1,\, w_2, \, w_3]$, $\sum_{i=1}^3 w_i=1$, $0 \leq w_i \leq 1 \, \forall i$. Assuming a single data point is chosen from and observed among the candidates $\Theta$ at each step, i.e. $\tilde x= \Omega_1$, a specific weighted sum formulation to address Problem~\ref{prob:search2} is obtained.

\begin{prop} \label{prob:weight1}
The solution, $\tilde x \in \Theta$, to the optimization problem
\begin{equation} \label{e:pw1}
 \max_{\tilde x \in \Theta} F(\tilde x)=\sum_{i=1}^3 F_i(\tilde x)= w_1 \hat f(\tilde x) 
- w_2 \dfrac{1}{N} \sum_{\tau \in \Theta} \abs{\hat e(\tilde x,\mc D, \tau)} 
+ w_3 \mc I (\tilde x, \hat f),
\end{equation}
constitutes the best search strategy for this weighted sum formulation of Problem~\ref{prob:search2}.
\end{prop}

As discussed in Subsection~\ref{sec:obsinfo} and stated in Proposition~\ref{thm:GPsearch2}, the information objective, $F_3$, in (\ref{e:pw1}) can be approximated by substituting it with GP variance $v(x)$ in (\ref{e:gp1}) to decrease computational load. Thus, an approximation to the solution in Proposition~\ref{prob:weight1} is:
\begin{prop} \label{prob:weight2}
The solution, $\tilde x \in \Theta$, to the optimization problem
\begin{equation} \label{e:pw2}
 \max_{\tilde x \in \Theta} F(x)= \sum_{i=1}^3 F_i(\tilde x)= w_1 \hat f(x) 
- w_2 \dfrac{1}{N} \sum_{\tau \in \Theta} \abs{\hat e(\tilde x,\mc D, \tau)} 
+ w_3 v(\tilde x),
\end{equation}
where $v(\tilde x)$ is defined in (\ref{e:gp1}), approximates the search strategy in Proposition~\ref{prob:weight1}.
\end{prop}

The weighting scheme described is only meaningful if the three objectives are of the same order of magnitude. Therefore, the original objective functions, $F_i$, $i=1,2,3$, have to be transformed or ``normalized''. There are many different approaches to perform such a transformation \cite{moosurvey1,moosurvey2}. The most common one, which coincidentally is known as normalization, aims to map each objective function to a predefined interval, e.g. $[0,\, 1]$. To do this, estimate first an upper $F_i^U$ and lower $F_i^L$ bound on each individual objective $F_i(x)$. Then, the $i^{th}$ normalized objective is
$$ F_i^N(x)= \dfrac{F_i(x) - F_i^L}{F_i^U -F_i^L }.$$

The main issue in normalization is to determine the appropriate upper and lower bounds, which is a very problem-dependent one. In the case of Proposition~\ref{prob:weight2}, the estimated functions $\hat f$ and $\hat e$ on the set $\Theta$ as well as the existing observations $\mc D$, can be utilized to obtain these values. The specific bounds for the respective objectives $F_1^U=\max_{x \in \Theta} \hat f(x)$, $F_1^L=\min_{x \in \Theta} \hat f(x)$, $F_2^U=\max_{x \in \Theta} \abs{\hat e(x,\mc D)}$, $F_2^L=0$, $F_3^U=\max_{x \in \Theta} \kappa(x)$, and $F_3^U=0$ provide a suitable starting estimate and can be combined with a prior domain knowledge if necessary. Thus, a normalized version of the formulation in Proposition~\ref{prob:weight2} is obtained.

\begin{prop} \label{prob:weight3}
The solution, $\tilde x \in \Theta$, to the optimization problem
\begin{equation} \label{e:pw3}
 \max_{\tilde x \in \Theta} F(x)= \sum_{i=1}^3 F_i^N(\tilde x)=  
  \dfrac{w_1}{\Delta_1} \left( \hat f(x) -  F_1^L\right) 
- \dfrac{w_2}{\Delta_2} \dfrac{1}{N} \sum_{\tau \in \Theta} \abs{\hat e(\tilde x,\mc D, \tau)} 
+ \dfrac{w_3}{\Delta_3} v(\tilde x),
\end{equation}
where $\Delta_i=F_i^U- F_i^L \; i=1,2,3$, provides an approximation to the best search strategy for solving the normalized weighted-sum formulation of Problem~\ref{prob:search2}.
\end{prop}


The \textbf{bounded objective function} method provides a suitable alternative to the weighted sum formulation above in addressing the multi-objective problem defined. The bounded objective function method minimizes the single most important objective, in this case $F_1(x)$, while the other two objective functions $F_2(x)$ and $F_3(x)$ are converted to  form additional constraints. Such constraints are in a sense similar to QoS ones that naturally exist in many real life problems  \cite{tcomm,alpcan-winet2,srikantbook}.
As an advantage, in the bounded objective formulation there is no need for normalization.

The bounded objective counterpart of the result in Proposition~\ref{prob:weight2} is as follows.
\begin{prop} \label{prob:bound1}
The solution, $\tilde x \in \Theta$, to the constrained optimization problem
\begin{eqnarray} \label{e:pw4}
 \max_{\tilde x \in \Theta} \hat f(x) \\
 \nonumber \text{such that } 0 \leq F_2(\tilde x)= \dfrac{1}{N} \sum_{\tau \in \Theta} \abs{\hat e(\tilde x,\mc D, \tau)} \leq b_1, \\
  \nonumber \text{and }  0 \leq F_3(\tilde x)= v(\tilde x) \leq b_2, 
\end{eqnarray}
where $b_1$ and $b_2$ are given (predetermined) scalar bounds on $F_2$ and $F_3$, respectively, provides an approximate best search strategy for a bounded-objective formulation of  Problem~\ref{prob:search2}.
\end{prop}

The advantage of the bounded objective function method is that it provides a bound on the information collection and estimation objectives while maximizing the estimated function. This leads in practice to an initial emphasis on information collection and correct estimation of the objective function. In that sense, the method is more ``classical'', i.e. follows the common method of learn first and maximize later. Furthermore, it does not require normalization, i.e. it is easier to deploy. The method has, however, a significant disadvantage which makes its usage
prohibitive. In large-scale or high-dimensional problems, the space to explore to satisfy any bound on information
is simply immense. Therefore, one does not have the luxury of identifying the function first to maximize it later as it would take too many samples to do this. In such cases, it makes more sense to deploy the weighted sum method, possibly along with a cooling scheme to modify the weights as part of a cooling scheme to balance depth-first vs.
breadth-first search.

Until now, it has been (implicitly) assumed that the static optimization problem at hand is stationary. However, in a variety of problems this is not the case and the function $f(x,t)$ changes with time. The decision making framework allows for modeling such systems in the following way. Let 
$$\mc O(t)=\{f(x_i,t_i)+n(x_i,t_i): x_i \in \mc D, t_i \leq t, \, \forall i\} ,$$
be the set of noisy or unreliable past observations until time $t$, where $n(x,t) \sim \mc N(0,\sigma(t))$ is the zero mean Gaussian ``noise'' term at time $t$. Now, the deterioration in the past information due to change in $f(x,t)$ can be captured by increasing the variance of the noise term, $\sigma(t)$, with time. For example, a simple linear dynamic can be defined as
$$ \dfrac{d\sigma(t)}{dt}= \eta,$$
where $\eta>0$ captures the level of stationarity, e.g. a large $\eta$ indicates a rapidly changing system and function $f(x,t)$.

\subsection{Algorithm}

An algorithmic summary of the solution approaches discussed above for a specific set of choices is provided by Algorithm~\ref{alg:algopt1}, which describes both weighted-sum and bounded objective variants. 
\begin{algorithm}[htbp]
   \caption{Optimization with Limited Information}
   \label{alg:algopt1}
\begin{algorithmic}[1]
  \STATE {\bfseries Input:} Function domain, $\mc X$, GP meta-parameters, objective weights $[w_1, w_2, w_3]$ or bounds $b_1, b_2$, initial data set $(\mc D, y)$.
  \STATE Use GP with a Gaussian kernel and specific expected error variances for function $\hat f$ and error function $\hat e$ estimation.
  \WHILE{Search budget available, $1 \leq n \leq N_{max}$.}
  \STATE Sample domain  $\mc X$ to obtain $\Theta(n)$. In some cases, $\Theta(n)=\Theta \; \forall n$.
  \STATE Estimate $\hat f$ and $\hat e$ based on observed data $(\mc D, y)$ on $\Theta(n)$ using GPs.
  \STATE Compute variance, $v(x)$, of $\hat f$ (\ref{e:gp1}) on $\Theta(n)$ as an estimate of $\mc I(\hat f)$.
    \IF{Weighted-sum method}
      \STATE  Next action maximizes a normalized and weighted sum of objectives $\sum_{i=1}^3 F_i^N$ as stated in Proposition \ref{prob:weight3}. 
    \ELSIF{Bounded objective method}
      \STATE  Next action is solution to the constrained problem in Proposition \ref{prob:bound1}.
    \ENDIF
  \STATE Update the observed data $(\mc D, y)$.
  \ENDWHILE
\end{algorithmic}
\end{algorithm}


\subsection{Numerical Analysis} \label{sec:numeric}

The Algorithm \ref{alg:algopt1} is illustrated next with multiple numerical examples. It is worth reminding that
the main issue here is to solve the optimization problems with minimum data using active learning. In all examples,
a uniform grid is used to sample the solution space rather than resorting to a more sophisticated method since
the examples are chosen to be only one or two dimensional for visualization purposes. 

\subsubsection*{Example 1}

The first numerical example aims to visualize the presented framework and algorithm. Hence, the chosen function is only one dimensional, $f(x)=sin(5x)/x$ on the interval $\mc X =[0.1, 3.9]$. The interval is linearly sampled to obtain a grid with a distance of $0.01$ between points, i.e. $\Theta=\{x_i \in \mc X \, \forall i: x_1=0.1, x_2=0.11,\ldots, x_N=3.9 \}$. A Gaussian kernel with variance $0.1$ is chosen for estimating both $\hat f$ and $\hat e$. The weights are equal to one, $w=[1,\, 1,\, 1]$, in the weighted-sum method. The bounds are $b_1=0.5$ for the error bound and $b_2=0.2$ for the bound on maximum variance estimate in the bounded objective method. The initial data consists of a single point, $x=0.1$.

\begin{figure}[htbp]
 \centering
 \includegraphics[width=0.8\columnwidth]{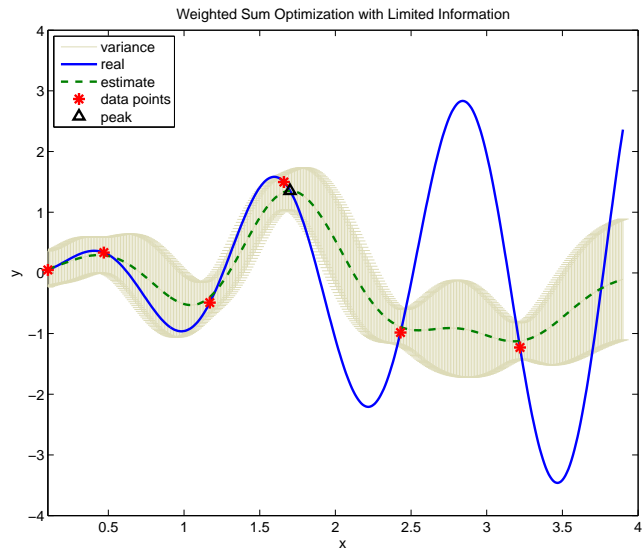}
 \caption{Optimization result using the weighted-sum method with $6$ data points.}
 \label{fig:weighted6}
\end{figure}

Figure~\ref{fig:weighted6} shows the results based on the normalized weighted-sum method in Proposition~\ref{prob:weight3} after $5$ iterations ($6$ samples in total, together with the initial data point). The variance here is $v(x)$ of the estimated function $\hat f$ using data points $\mc D$. Clearly, the estimated peak is not the one of the real function $f$.

\begin{figure}[htbp]
 \centering
 \includegraphics[width=0.8\columnwidth]{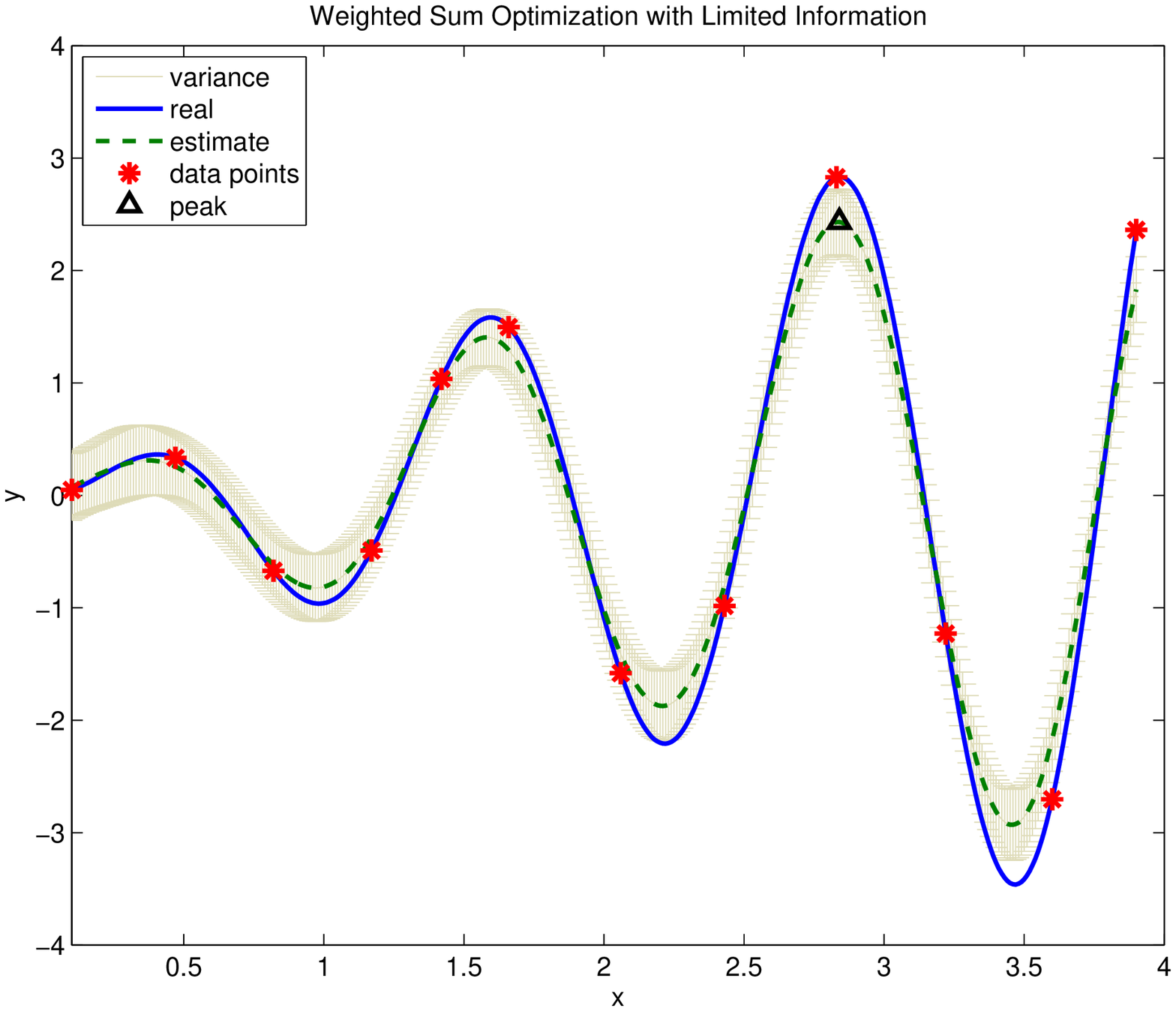}
 \caption{Optimization result using the weighted-sum method with $12$ data points.}
 \label{fig:weighted12}
\end{figure}

Next, Figure~\ref{fig:weighted12} shows that after $11$ iterations ($12$ data points in $\mc D$), the function and the location of  its peak is estimated correctly. 
The sequence of points selected during the iteration process are: 
$$ \mc D=\{0.47, 3.22, 1.17, 1.66, 2.43, 2.06, 3.9, 2.83, 3.6, 0.82, 1.42 \}.$$

\begin{figure}[htbp]
 \centering
 \includegraphics[width=0.8\columnwidth]{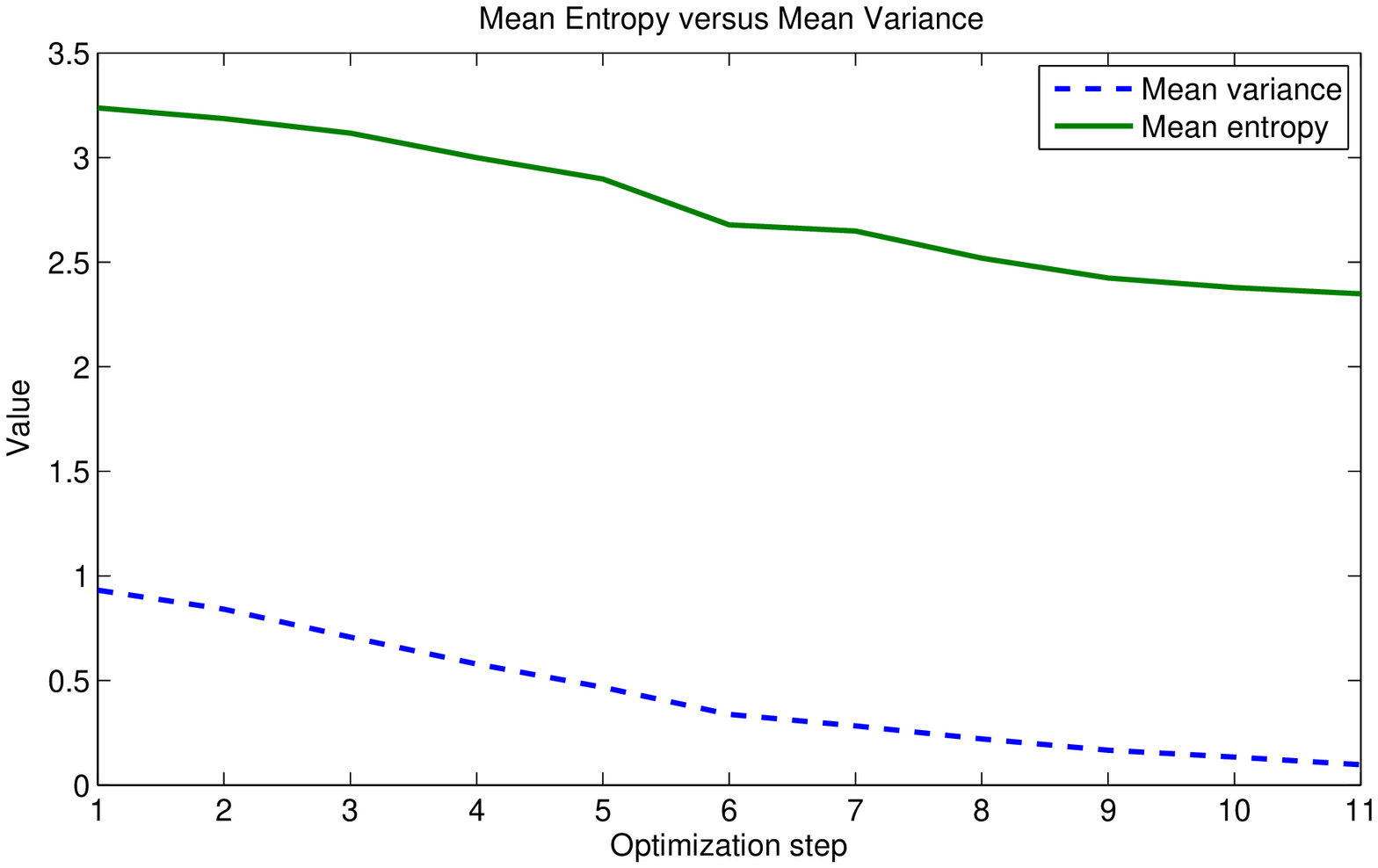}
 \caption{Mean variance $v$ and entropy $\mc I$ on $\Theta$ at each iteration step.}
 \label{fig:info1}
\end{figure}
\begin{figure}[htbp]
 \centering
 \includegraphics[width=0.8\columnwidth]{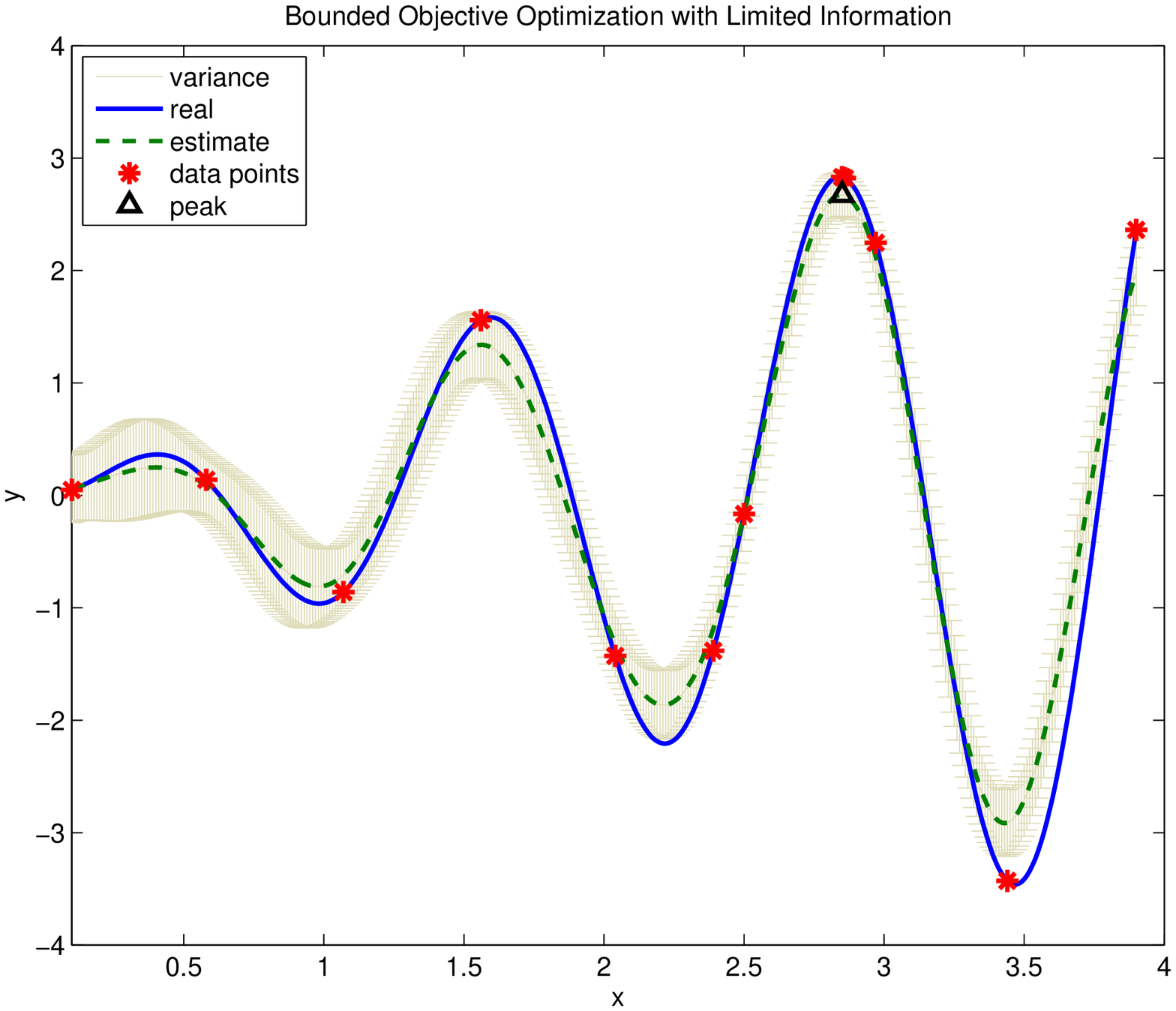}
 \caption{Optimization result using the bounded objective method with $12$ data points.}
 \label{fig:bounded1}
\end{figure}

The amount of information obtained during the iterative optimization is of particular interest. Figure~\ref{fig:info1} depicts the mean variance $v$ and entropy $\mc I$ of the estimated function $\hat f$ on  $\Theta$ at each iteration step. In this specific example, the two quantities are very well correlated. Note, however, that this correlation is a function of the relative weights between information collection and other objectives.

Finally, Figure~\ref{fig:bounded1} depicts the results of the bounded objective method with the given bounds. The number of iterations is $11$ as before, which gives an opportunity of direct comparison with the weighted-sum method.
The sequence of points selected during the iteration process are: 
$$ \mc D=\{0.47, 3.22, 1.17, 1.66, 2.43, 2.06, 3.9, 2.83, 3.6, 0.82, 1.42 \}.$$

\subsubsection*{Example 2}

\begin{figure}[htbp]
 \centering
 \includegraphics[width=0.8\columnwidth]{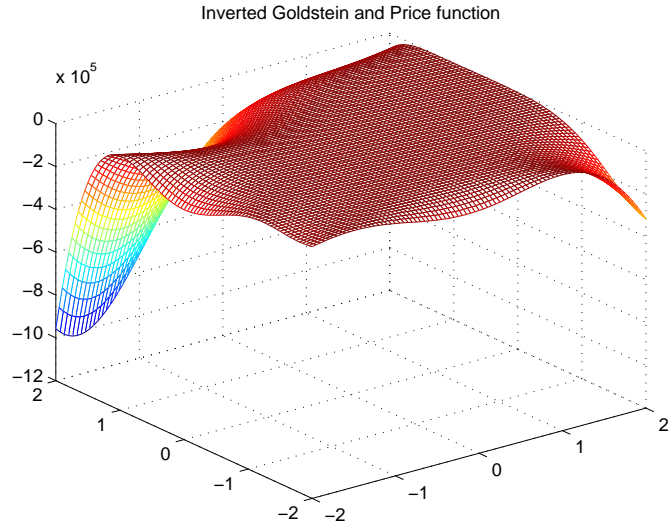}
 \caption{The inverted Goldstein\&Price function~\cite{gpfunc}.}
 \label{fig:gpfunc}
\end{figure}

\begin{figure}[htbp]
 \centering
 \includegraphics[width=0.8\columnwidth]{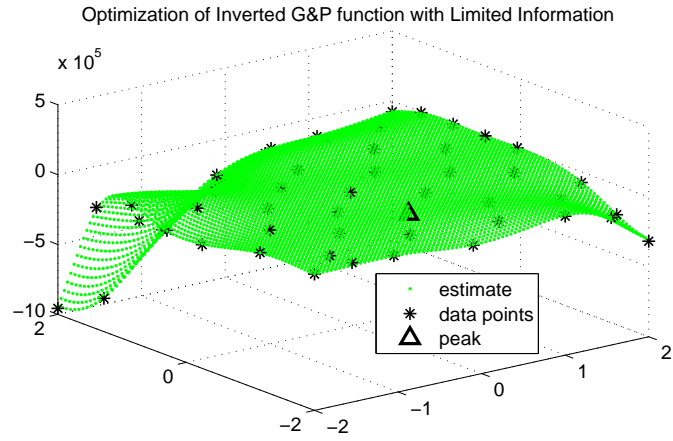}
 \caption{Optimization of the inverted Goldstein\&Price function~\cite{gpfunc} using the weighted-sum method with $50$ data points.}
 \label{fig:gpopt1}
\end{figure}

\begin{figure}[htbp]
 \centering
 \includegraphics[width=0.8\columnwidth]{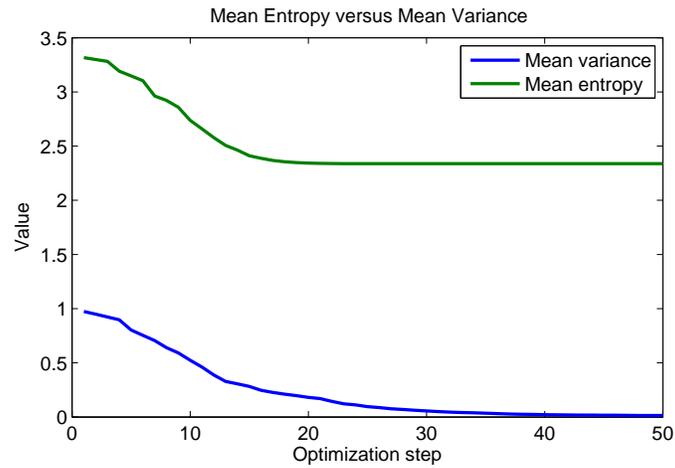}
 \caption{Mean variance $v$ and entropy $\mc I$ on $\Theta$ at each iteration step.}
 \label{fig:infogp}
\end{figure}

The objective function in the second numerical example is the  Goldstein\&Price function~\cite{gpfunc}, which
is shown in Figure~\ref{fig:gpfunc} in its inverted form to ensure consistency with the maximization formulation in this paper.
The problem domain consists of the two dimensional rectangular region $\mc X=[-2, 2] \times [-2, 2]$, which is linearly sampled to obtain a uniform grid with a $0.05$ interval between sample points. A Gaussian kernel with variance $0.5$ and $0.1$ is chosen for estimating $\hat f$ and $\hat e$, respectively. The weighted-sum method is utilized in Algorithm~\ref{alg:algopt1} with the weights $w=[4,\, 2,\, 3]$. The search budget is chosen as $50$ before stopping the algorithm (for the search space of approx. $6400$ samples in the grid). The real global minimum (peak) of the (inverted) Goldstein\&Price function is at $(0, -1)$ and the location found by the algorithm using the $50$ data points is $(-0.15, -1.05)$. Figure~\ref{fig:gpopt1} depicts the estimated function,
the data points as well as the optimum found. Although the real optimum value is $-3$ (in the inverted version) while the obtained one is $-9.75$, the result is still very satisfactory considering that the simple sampling scheme used and the Goldstein\&Price function takes values in a range of $1$ million, i.e. the error is less than $0.001$ percent of the range. Finally, Figure~\ref{fig:infogp} depicts the mean variance $v$ and entropy $\mc I$ of the estimated function  $\hat f$  on  $\Theta$ at each iteration step.

\subsubsection*{Example 3}

\begin{figure}[htbp]
 \centering
 \includegraphics[width=0.8\columnwidth]{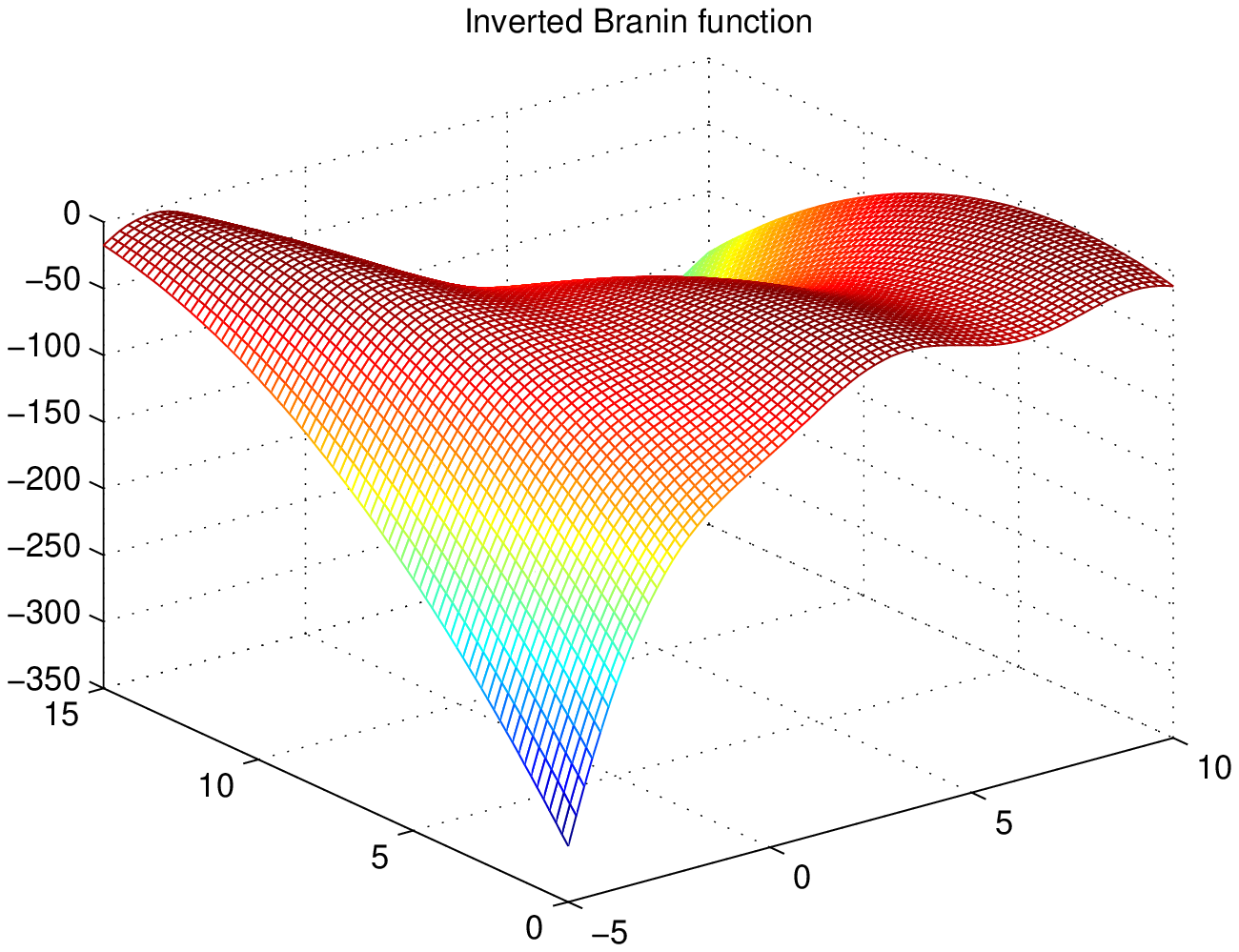}
 \caption{The inverted Branin function~\cite{braninfunc}.}
 \label{fig:braninfunc}
\end{figure}

The third example uses the same setup as the second one but this time with the (inverted) Brain function~\cite{braninfunc}
shown in Figure~\ref{fig:braninfunc}. The rectangular problem domain $\mc X=[-5, 10] \times [0, 15]$ is sampled uniformly to obtain a grid of points with a $0.2$ interval. The real global minimums (peaks) of the (inverted) Branin function are at 
$(9.4,2.47)$, $(-\pi,12.28)$, and $(\pi,2.28)$ whereas the locations found by the algorithm are $(9,2.6)$, $(-3.2,12)$, and $(3,2.2)$. The values at these locations found vary between $-4.3$ and $-0.5$ compared to the real global value of
$-0.4$ (of the inverted function). Thus, the algorithm again performs satisfactorily. Figure~\ref{fig:gpopt1} shows the
computed location of one optimum, the data points, as well as the estimated function based on the data points.
\begin{figure}[htbp]
 \centering
 \includegraphics[width=0.8\columnwidth]{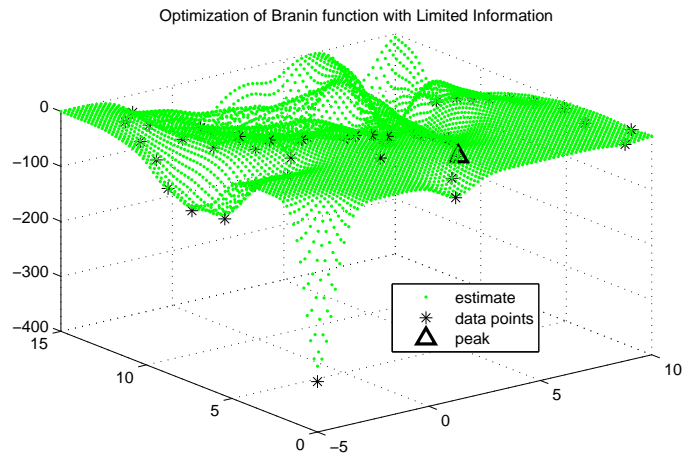}
 \caption{Optimization of the inverted Branin function~\cite{braninfunc} using the weighted-sum method with $50$ data points.}
 \label{fig:braninopt1}
\end{figure}

\subsubsection*{Example 4}

\begin{figure}[htbp]
 \centering
 \includegraphics[width=0.8\columnwidth]{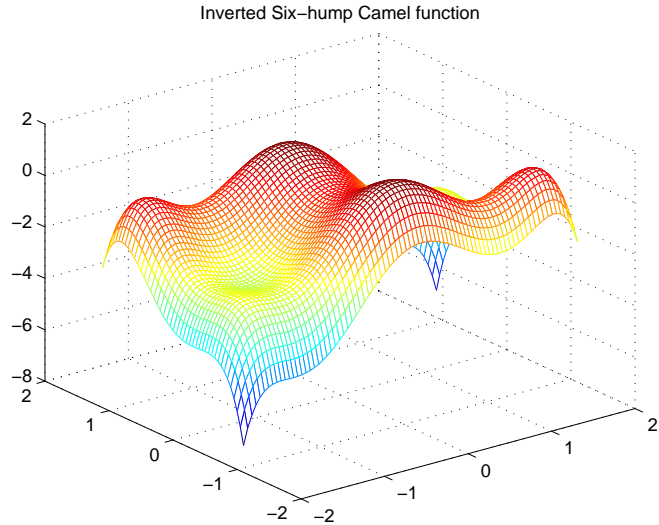}
 \caption{The inverted six-hump camel function.}
 \label{fig:camfunc}
\end{figure}

\begin{figure}[htbp]
 \centering
 \includegraphics[width=0.8\columnwidth]{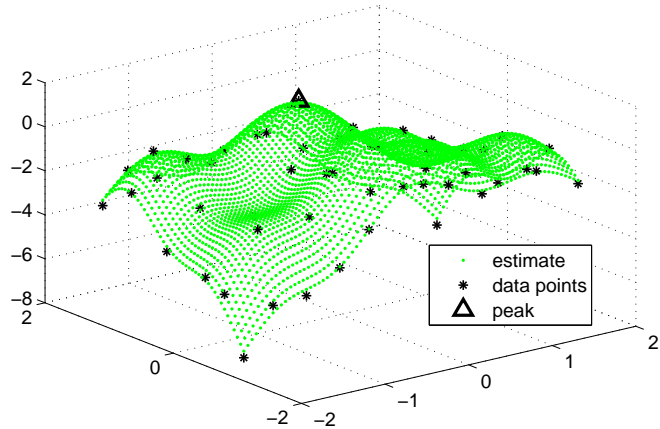}
 \caption{Optimization of the inverted six-hump camel function~\cite{globaloptold} using the weighted-sum method with $50$ data points.}
 \label{fig:camopt1}
\end{figure}

The fourth example is based on the six-hump camel function~\cite{globaloptold} (see Figure~\ref{fig:camfunc}) on the domain $\mc X=[-2, 2] \times [-2,2]$, which is sampled uniformly with a $0.05$ interval. All of the parameters are chosen to be the same as before. Figure~\ref{fig:camopt1} shows the
computed location of two optimums, the $50$ data points, as well as the estimated function based on the data points.
The optimum locations found are $(0, \, 0.65)$ and $(0.05,\, -0.6)$ with respective values of $0.98$ and $1.06$,
whereas the real locations are $(-0.09,0.71)$ and $(0.09,-0.71)$ with the value $1.03$.

\section{Literature Review} \label{sec:literature}

Decision making with limited information is related to search theory. The idea of using information (theory) in this context is hardly new as evidenced by the article ``A New Look at the Relation Between Information Theory and Search Theory'' from 1979 \cite{pierce}. The subject is further studied in \cite{jaynes}. The topic of optimal search is more recently revisited by \cite{Zhu-search}, which contains substantial historical notes and studies problems where the search target distribution in itself is unobservable.

The book \cite{MacKaybook} provides important and valuable insights into the relationship between information theory, inference, and learning. Measuring information content of experiments using Shannon information is explicitly mentioned and a slightly informal version of the bisection example in Subsection~\ref{sec:obsinfo} is discussed. However, focusing mainly on more traditional coding, communication, and machine learning topics, the book does not discuss the type of decision making problems presented in this paper. 

Learning plays an important role in the presented framework, especially \textit{regression}, which is a classical machine (or statistical) learning method. A very good introduction to the subject can be found in \cite{Bishopbook}. A complementary and detailed discussion on kernel methods is in \cite{schoelkopfbook}. Another relevant topic is Bayesian inference \cite{Tipping,MacKaybook}, which is in the foundation of the presented framework. In machine learning literature, Gaussian processes (GPs) are getting increasingly popular due to their various favorable characteristics. The book \cite{GPbook} presents a comprehensive treatment of GPs. Additional relevant works on the subject include \cite{MacKaybook,schoelkopfbook,MacKayGP}, which also discuss GP regression.

Convex optimization \cite{boydbook} is a well-understood topic that is often easy to handle even if available information is limited. Optimizing nonconvex functions, however, is still a research subject \cite{globaloptsurvey}. It is interesting to note that the method known as \textit{kriging} in global optimization is almost the same as GP regression in machine learning. The field \textit{stochastic programming} focuses on optimization under uncertainty but assumes a certain amount of prior knowledge on the problem at hand and models the uncertainty probabilistically \cite{stochasticsurvery}. The popular heuristic method \textit{simulated annealing} \cite{simannealing} is essentially based on iterative random search. Another popular heuristic scheme particle swarm optimization \cite{swarmorig} is also based on random search but parallel in nature as a distinguishing characteristic rather than iterative. 

Gaussian processes have been recently applied to the area of optimization and regression \cite{BoylePhD} as well as system identification \cite{ThompsonPhD}. While the latter mentions active learning, neither work discusses explicit information quantification or builds a connection with Shannon information theory. The recent articles \cite{GPopt,turner1}, which utilize GP regression for optimization  in a setting similar to the one in this paper and for state-space inference and learning, respectively, do not consider information-theoretic aspects of the problem, either. Likewise,
the article \cite{kriging1} on stochastic black box optimization, which considers a problem similar to the one here, does not take into account explicit measurement of information.

The area of active learning or experiment design focuses on data scarcity in machine learning and makes use of Shannon information theory among other criteria \cite{activelearning}. The paper \cite{MacKaydataselect} discusses objective functions which measure the expected informativeness of candidate measurements within a Bayesian learning framework. 
The subsequent study \cite{gpactive1} investigates active learning for GP regression using variance as a (heuristic) confidence measure  for test point rejection.


\section{Discussion} \label{sec:discuss}

The foundation of the approach adopted in this paper is Bayesian inference, where the main idea is to choose an a priori model and update it with actual experimental data observed (see \cite[Chap. 2]{MacKaybook} for a beautiful introductory discussion on the subject). As long as the a priori model is close to the reality (of the problem at hand), this inference methodology works very efficiently as indicated by  the numerical examples in Section~\ref{sec:numeric}. In many cases this background information, which is sometimes referred to as ``domain knowledge'', is already available. However, in others one has to explore the model domain and learn model meta-parameters in a time scale naturally longer than the one of actual optimization \cite{MacKayGP}. 

The GP regression adopted in the presented framework is only one method for function estimation and other, e.g. parametric, methods can easily replace GP for the regression part. In any case, the regression methodology here is consistent with
the principle of ``Occam's razor'', more specifically its interpretation using Kolmogorov complexity \cite{algobook}. A priori,
the optimization problems at hand 
are more probable to be simple
rather than complex to describe in accordance with \textit{universal distribution} \cite{algobook}. Hence, given a data set it is reasonable to start describing it with the simplest explanation. GP regression already incorporates this line of thinking by relying on a kernel-based approach and making use of the representer theorem \cite[Chap. 6.2]{GPbook}. As a visual example, we refer Figures \ref{fig:weighted6} and \ref{fig:weighted12} for a comparison of function estimates with different sets of available data.

This paper considers a class of problems where data is scarce and obtaining it is costly. Information theory plays an especially important role in devising optimal schemes for obtaining new data points (active learning). The entropy measure from Shannon information theory provides the necessary metric for this purpose, which quantifies the ``exploration'' aspect of the problem. Using a multi-objective optimization formulation, the presented framework allows explicit weighting of \textit{exploration} vs. \textit{exploitation} aspects. This trade-off is also very similar to one between the well-known depth-first vs. breadth-first search algorithms in search theory.

The amount of information obtained from each data point is different here only because a specific a priori general model is utilized to explain the observed data (GP regression). Because of this the amount of information obtained is specific to the model. Otherwise, without this Bayesian approach, each data point would give the same information (inversely proportional to the total number of candidate points).

The illustrative examples discussed are low-dimensional, which makes it possible to use grids for sampling. However, in higher dimensions (i.e. when the problem is much more ``difficult'') this ``luxury'' is not affordable and one has to necessarily resort to Monte Carlo methods. In such cases, the trade-off between exploration and exploitation
is even more emphasized. Possible methods to address this issue include, ``cooling'' approaches similar
to those used in simulated annealing, multi-resolution sampling based on region of interest or using topological properties of Gaussian mixtures to intelligently estimate candidate points based on the current state.

The optimization approach presented here can also be interpreted from a biological perspective. If an analogy between
the decision-maker and a biological organism is established, then the a-priori Bayesian model (meta parameters of the GP) that is refined over a long time scale corresponds to  evolution of a species in an environment
(problem domain). Each individual organism belogning to the species obtains new information to achieve its objective while preserving resources as much as possible. The existing evolutionary basis (GP model) gives them an advantage to find a solution much faster compared to random search. From the perspective of the species, it also makes sense for some of its members to explore the model (meta parameter) domain and further refine it through adaptation. Those with better meta parameters achieve then their objectives even more efficiently and obtain an evolutionary edge in natural selection (assuming competition).

\section{Conclusion} \label{sec:conclusion}

The decision making framework presented in this paper addresses the problem of decision making under limited information by taking into account the information collection (observation), estimation (regression), and (multi-objective) optimization aspects in a holistic and structured manner. The methodology is based on Gaussian processes
and active learning. Various issues such as quantifying information content of new data points using information theory, the relationship between information and GP variance as well as related approximation and multi-objective optimization schemes are discussed. The framework is demonstrated with multiple numerical examples.

The presented framework should be considered mainly as an initial step. Future research directions are abundant
and include further investigation of the exploration-exploitation trade-off, adaptive weighting parameters, and random sampling methods for problems in higher dimensional spaces. Additional research topics are the relationship of the framework with genetic/evolutionary methods, dynamic control problems, and multi-person decision making, i.e.
game theory.

%
%

\section*{Acknowledgements} 
This work is supported by Deutsche Telekom Laboratories. The author wishes to thank Lacra Pavel, Slawomir Stanczak, Holger Boche, and Kivanc Mihcak for stimulating discussions on the subject.



\end{document}